\documentclass[12pt]{amsart}
\usepackage{amsthm,amsfonts, amssymb, amscd}
\newtheorem{definition}{Definition}[section]
\newtheorem{thm}{Theorem}[section]
\newtheorem{prop}{Proposition}[section]

\newtheorem{cor}{Corollary}[section]
\def\endproof{$\square$}
\def\Mg{{\mathcal M}_{g}}
\def\Mgbar{\overline{{\mathcal M}}_g}

\def\C{{\mathbb C}}

\def\S{{\mathcal S}}
\def\O{{\mathcal O}}

\def\P{{\bf P}}
\def\A{{\bf A}}

\def\Hbar{\overline{H}}
\def\PGL{{\bf P}GL}
\def\Z{{\mathbb Z}}

\begin{document}
\title[Notes on moduli]
{Notes on the construction of the moduli space of curves}
\author{Dan Edidin}
\thanks{Latest Revision - May 1998}
\address{Department of Mathematics\\ University of Missouri\\
Columbia MO 65211}
\maketitle
The purpose of these notes is to discuss the problem of
moduli for curves of genus
$g \geq 3$
\footnote{Because every curve of genus 1 and 2 has non-trivial
automorphisms, the problem of moduli is more subtle in this case than
for curves of higher genus} and
outline the construction of the (coarse) moduli scheme of
stable curves due to Gieseker. We present few complete proofs.
Instead we try and explain the subtleties and give precise references
to the literature. The notes are broken into 4 parts.

In Section 1 we discuss the general problem of constructing
a moduli ``space'' of curves. We will also state
results about its properties, some of which will be discussed
in the sequel.

We begin Section 2 by recalling from \cite{DM} (and also \cite{Vi}) the
definition of a groupoid, and
define the moduli groupoid of curves, as well as the quotient groupoid
of a scheme by a group.
We then give the conditions required for a groupoid
to be a stack, and prove that the quotient groupoid
of a scheme by a group is a stack. After discussing properties
of morphisms of stacks, we define a Deligne-Mumford stack
and prove that if a group acts on a scheme so that
the stabilizers of geometric points are finite and reduced
then the quotient stack is Deligne-Mumford. We conclude
Section 2 with the definition of the moduli space
of a Deligne-Mumford stack.

In Section 3 the notion of a stable curve is introduced, and we
define the groupoid of stable curves. The groupoid of smooth curves
is a sub-groupoid.
We then prove that the groupoid  of stable curves of genus $g\geq 3$
is equivalent to the quotient groupoid of a
Hilbert scheme by the action of the projective linear
group with finite, reduced stabilizers at geometric points.
By the results of the previous section we can
conclude the the groupoid of stable curves is a Deligne-Mumford stack
defined over $Spec\; \Z$
(as is the groupoid of smooth curves). We also
discuss the results of \cite{DM} on the irreducibility
of the moduli stack.

In Section 4 we prove that a geometric quotient of
a scheme by a group is the moduli space for the quotient
stack.
We then discuss the method of geometric invariant theory
for constructing geometric quotients for the actions of
reductive groups.
Finally, we briefly outline Gieseker's approach to
constructing the coarse moduli space over an algebraically closed
field
as the quotient of the aforementioned Hilbert scheme.

{\em Acknowledgments:} These notes are based on lectures the author
gave at the Weizmann Institute, Rehovot, Israel in July 1994.  It is a
pleasure to thank Amnon Yekutieli and Victor Vinnikov for the
invitation, and for many discussions on the material in these notes.
Thanks also to Alessio Corti and Andrew Kresch for useful discussions.
Comments and corrections are welcome!

\section{Basics}
\begin{definition}
Let $S$ be a scheme. By a smooth curve of genus $g$ over $S$ we mean
a proper, flat, family $C \rightarrow S$ whose geometric
fibers are smooth, connected 1-dimensional schemes of genus $g$.
\end{definition}

Note: By the genus of a smooth, connected curve $C$ over
an algebraically closed field, we mean $\mbox{dim }H^0(C,\omega_C)
= \mbox{dim }H^1(C,{\mathcal O}_C)=g$, where $\omega_C$ is the sheaf
of regular 1-forms on $C$.
If the ground field is $\C$, then $C$ is a smooth, compact Riemann
surface, and the algebraic definition of genus is the same as the
topological one.

The basic problem of moduli is to classify curves of genus $g$.
As a start, it is desirable to construct a space $\Mg$ whose
geometric points represent all possible isomorphism classes
of smooth curves. In the language of complex varieties, we are looking
for a space that parametrizes all possible complex structures
we can put on a fixed surface of genus $g$.

However, as modern (post-Grothendieck) algebraic geometers
we would like $\Mg$ to have further functorial properties.
In particular given a scheme $S$, a curve $C \rightarrow S$
should correspond to a morphism of $S$ to $\Mg$
(when $S$ is the spectrum of an algebraically closed field
this is exactly the condition
of the previous paragraph).

In the language of functors, we are trying to find
a scheme $\Mg$ which represents the functor
${\mathcal F}_{\Mg}: \mbox{Schemes} \rightarrow \mbox{Sets}$
which assigns to a scheme $S$ the set of isomorphism
classes of smooth curves of genus $g$ over $S$.

Unfortunately, such a moduli space can not exist
because some curves have non-trivial
automorphisms. As a result,  it is possible to construct non-trivial
families $C \rightarrow B$ where each fiber has the same
isomorphism class. Since the image of $B$ under the corresponding
map to the moduli space is a point, if the moduli space represented
the functor ${\mathcal F}_{\Mg}$ then $C \rightarrow B$ would be
isomorphic to the trivial product family.

Given a curve $X$ and a non-trivial (finite) group $G$ of automorphisms
of $X$ we construct a non-constant family $C \rightarrow B$
where each fiber is isomorphic to $X$ as follows. Let $B'$ be
a scheme with a free $G$ action, and let $B = B'/G$ be the quotient.
Let $C'=B' \times X$. Then $G$ acts on $C'$ by acting as
it does on $B'$ on the first factor, and by automorphism on the
second factor. The quotient $C'/G$ is a family of curves over
$B$. Each fiber is still isomorphic to $X$, but $C$ will not
in general be isomorphic to $X \times B$.

There is however, a {\em coarse} moduli scheme
of smooth curves. By this we mean:
\begin{definition} (\cite[Definition 5.6, p.99]{GIT})
There is a scheme $\Mg$ and a natural transformation
of functors $\phi:{\mathcal F}_{\Mg} \rightarrow Hom(\_\_,\Mg)$
such that

1. For any algebraically closed field $k$,
the map $\phi:{\mathcal F}_{Mg}(\Omega) \rightarrow Hom(\Omega,\Mg)$
is a bijection, where $\Omega= Spec \; k$.

2. Given any scheme $M$ and a transformation $\psi:{\mathcal F}_{\Mg}
\rightarrow Hom(\_\_,M)$ there is a unique transformation
$\chi:Hom(\_\_,\Mg) \rightarrow Hom(\_\_,M)$ such that $\psi=\chi \circ \phi$.

\end{definition}
The existence of $\phi$ means that given a family of
curves $C \rightarrow B$ there is an induced map to $\Mg$.
We do not require however, that a map to moduli gives a
family of curves (as we have already seen a non-constant
family with iso-trivial fibers).
However, condition 1 says that giving a curve over an algebraically
closed field is equivalent to giving a map of that field into
$\Mg$.

Condition 2 is imposed so that the moduli space is a universal
object.

In his book on geometric invariant theory Mumford proved the following
theorem.
\begin{thm}(\cite[Chapter 5]{GIT})
Given an algebraically closed field $k$ there
is a coarse moduli scheme $\Mg$ of dimension $3g-3$
defined over $Spec \; k$, which is quasi-projective and irreducible.
\end{thm}

The proof of this theorem will be subsumed in our general
discussion of the construction of $\Mgbar$, the moduli
space of stable curves.

A natural question to ask at this point, is whether
$\Mg$ is complete (and thus projective). The answer
is no. It is quite easy to construct curves
$C \rightarrow Spec\; {\mathcal O}$ with $\O$ a D.V.R. which
has function field $K$, where
$C$ is smooth, but the fiber over the residue
field is singular. Since
$C$ is smooth, the restriction $C_K \rightarrow Spec\;K$
is a smooth curve over $Spec\;K$, so there
is a map $Spec\;K \rightarrow \Mg$.
The existence of such a family
does not prove anything, since we must show
that we can not replace the special fiber
by a smooth curve. The total space $\tilde C$ of
a family with modified special fiber is birational to
$C$. Since $C$ and $\tilde C$ are surfaces (being
curves over 1-dimensional rings) there must be a sequence
of birational transformations (centered in the special fiber)
taking one to the other.
It appears therefore, that it suffices to construct
a family $C \rightarrow Spec\;\O$ such that no birational
modification of $C$ centered in the special fiber
will make it smooth.
Unfortunately, because $\Mg$ is only a coarse moduli scheme,
the existence of such a family does not prove that $\Mg$ is
incomplete. The reason is that there may be a map
$Spec\;{\mathcal O} \rightarrow \Mg$ extending the original map $Spec\;K
\rightarrow \Mg$ without there being a family of smooth curves
$\tilde{C} \rightarrow Spec\;{\mathcal O}$ extending
$C_K \rightarrow Spec\; K$. However, we will see when we discuss
the valuative criterion of properness for Deligne-Mumford stacks that it
suffices to show that for every finite extension $K \subset K'$
we can not complete the induced family
of smooth curves $C_{K'} \rightarrow Spec\; K'$
to a family of smooth curves $C' \rightarrow Spec\;{\mathcal O'}$,
where ${\mathcal O'}$ is the integral closure of ${\mathcal O}$
in $K'$.

{\bf Example:} The following family shows that ${\mathcal M}_3$
is not complete. It can be easily generalized to
higher genera.
Consider the  family $x^4 +xyz^2 + y^4 +t(z^4 + z^3x + z^3y + z^2y^2)$
of quartics in $\P^2 \times Spec\;{\mathcal O}$ where ${\mathcal O}$
is a D.V.R. with uniformizing parameter $t$.
The total space
of this family is smooth, but over the closed point
the fiber is a quartic with a node at the point
$(0:0:1) \in \P^2$. Moreover, even after base change,
any blow-up centered at the singular point of the
special fiber contains a $(-2)$ curve, so there is no
modification that gives us a family of smooth curves.

Since $\Mg$ is not complete, a natural question
is to ask whether it is affine. The answer again
is no. This follows from the fact
that $\Mg$ has a projective compactification
in the
moduli of abelian varieties, such that boundary has
codimension 2. In particular, this means
that there are complete curves in $\Mg$. On
the other hand, Diaz proved the following theorem
(\cite{Di}).
\begin{thm}
Any complete subvariety of $\Mg$ has dimension
less than $g-1$.
\end{thm}
It is not known how close this bound is to being sharp.

Finally we state a spectacular theorem due largely to
Harris and Mumford.
\begin{thm}
For $g>23$, $\Mg$ is of general type.
\end{thm}
The importance of this theorem is that
until its proof, people believed that
$\Mg$ was rational, or at least unirational.
The reason for this belief was that
for $g \leq 10$, the (uni)rationality of $\Mg$
had been affirmed by the Italian school.

\section{Stacks}
Let $S$ be a scheme, and let $\S=(Sch/S)$ be the category
of schemes over $S$.
\subsection{Groupoids}
\begin{definition}
A category over $S$ is a category $F$ together with
a functor $p_F:F\rightarrow \S$. If $B \in Obj(\S)$
we say $X$ lies over $B$ if $p_F(X)=B$.
\end{definition}

\begin{definition} (see also \cite[Definition 7.1]{Vi})
If $(F,p_F)$ is a category over $S$, then it
is a groupoid over $S$ if the following conditions hold:

(1) If $f:B' \rightarrow B$ is a morphism in $\S$, and
$X$ is an object of $F$ lying over $B$, then
there is an object $X'$ over $B'$ and a morphism $\phi:X' \rightarrow X$
such that $p_F(\phi) = f$.

(2) Let $X,X',X''$ be objects of $F$ lying over
$B,B',B''$ respectively.
If $\phi:X' \rightarrow X$ and $\psi:X'' \rightarrow X$
are morphisms in $F$, and $h:B' \rightarrow B''$ is a
morphism such that $p_F(\psi) \cdot h = p_F(\phi)$ then
there is a unique morphism $\chi:X' \rightarrow X''$
such that $\psi\cdot \chi =\phi$ and $p_F(\chi) = h$.

\end{definition}

Note that condition (2) implies that a morphism $\phi:X' \rightarrow
X$ of objects over $B'$ and $B$ respectively is an isomorphism
if and only
if $p_F(\phi):B' \rightarrow B$ is an isomorphism. (To see that
$p_F(\phi)$ being an isomorphism is sufficient to ensure that
$\phi$ is an isomorphism, apply condition (2) where one of
the maps is $p_F(\phi)$ and the other the identity, and
lift $p_F(\phi)^{-1}:B \rightarrow B'$ to $\phi^{-1}:X \rightarrow
X'$. The other direction is trivial.)
Define $F(B)$ to be the subcategory
consisting of all objects $X$ such that $p_F(X) = B$ and morphisms
$f$ such $p_F(f) = id_B$. Then $F(B)$ is a groupoid; i.e. a category
where all morphisms are isomorphisms.
This is the reason we say $F$ is a groupoid over
$S$.
Also note that condition (2) implies that
the object $X'$ over $B'$ in condition (1) is unique up to canonical
isomorphism.
This object will be called the pull-back of $X$ via
$f$ and denoted $f^*X$.

{\bf Example:} If $F:\S \rightarrow \mbox{Sets}$ is a
contravariant functor, then we can associate a
groupoid (also called $F$) whose objects are
pairs $(B,\beta)$ where $B$ is an object of $\S$
and $\beta \in F(B)$. A morphism $(B',\beta^{\prime})
\rightarrow (B,\beta)$ is an $S$-morphism $f:B' \rightarrow B$
such that $F(f)(\beta) = \beta^{\prime}$. In this
case $F(B)$ in the groupoid sense is just the
set $F(B)$ in the functor sense; i.e. all morphisms
in the groupoid $F(B)$ are the identity.
In particular, any $S$-scheme $B$ is
a groupoid via its functor of points $Hom(\_,B)$.

{\bf Example:} If $X/S$ is a scheme and $G/S$
is a group scheme acting on $X$
then we define the quotient groupoid
$[X/G]$ as follows. The sections of $[X/G]$
over $B$ are $G$-principal bundles $E \rightarrow B$
together with a $G$-equivariant map $E \rightarrow X$.
A morphism from $E' \rightarrow B'$ to $E \rightarrow B$
is a commutative diagram
$$\begin{array}{ccc}
E' & \rightarrow & E\\
\downarrow & & \downarrow\\
B'&  \rightarrow & B
\end{array}$$
such that $E' \simeq E \times_B' B$.
If the action is free and a quotient scheme $X/G$ exists,
then
there is an equivalence of categories between
$[X/G]$ and the groupoid associated to the scheme
$X/G$.

{\bf Example:} Of particular importance to us is the groupoid
$F_{\Mg}$ defined over $Spec\;\Z$. The objects of
$F_{\Mg}$ are smooth curves as defined in part 1. A morphism
from $X' \rightarrow B'$ to $X \rightarrow B$ is
a commutative diagram
$$\begin{array}{ccc}
X' & \rightarrow & X\\
\downarrow & & \downarrow\\
B'&  \rightarrow & B
\end{array}$$
which induces an isomorphism $X' \simeq X \times_B B'$.
The functor $\Mg \rightarrow \mbox{Sch/}\Z$ sends
$X \rightarrow B$ to $B$. We will eventually prove that
$F_{\Mg}$ is a quotient groupoid as in the previous
example.

{\bf Warning:} The groupoid we have just defined
is not the groupoid associated to the moduli functor
we defined in Part 1. The groupoid here is not
a functor, since if $X/B$ is a curve with non-trivial
automorphisms $F(B)$ will not be a set because
there are morphisms which are not the identity. (A set is
a groupoid where all the morphisms are the identity.)

\subsection{Definition of a stack}

Let $(F,p_F)$ be a groupoid.
Let $B$ be an $S$-scheme and let $X$ and $Y$ be
any objects in $F(B)$. Define a functor
$Iso_B(X,Y):(\mbox{Sch}/B) \rightarrow (\mbox{Sets})$
by associating to any morphism $f:B' \rightarrow B$, the set of isomorphisms
in $F(B')$ between $f^*X$ and $f^*Y$.

If $X=Y$ then
$Iso_B(X,X)$ is the functor whose sections
over $B'$ mapping to $B$ are the automorphisms of
the pull-back of $X$ to $B'$.

In the case of curves Deligne and Mumford proved that $Iso_B(X,Y)$
is represented by a scheme $\mbox{{\bf Iso}}_B(X,Y)$,
because $X/B$ and $Y/B$
have canonical polarizations (\cite[p.84]{DM}).
When $X=Y$ then Deligne and Mumford prove directly
that the $\mbox{{\bf Iso}}_{B}(X,X)$ is finite and unramified
over $B$ (\cite[Theorem 1.11]{DM}). Applying the
theorem to  $B= Spec\;k$,
where $k$ is an algebraically closed field, this
theorem proves that every curve has a finite automorphism
group.

Note that the scheme $\mbox{{\bf Iso}}_B(X,X)$
need not be flat over $B$. For example,
if $X/B$ is a family of curves, the number of
points in the fibers of $\mbox{{\bf Iso}}_B(X,X)$ over $B$
will jump over the points $b \in B$ where the fiber
$X_b$ has non-trivial automorphisms.

\begin{definition}
A groupoid $(F,p_F)$ over $S$ is a stack if

(1) $Iso_B(X,Y)$ is a sheaf in the
\'etale topology for all $B$, $X$ and $Y$.

(2) If $\{B_i \rightarrow B\}$ is a covering of
$B$ in the \'etale topology, and $X_i$ is
a collection of objects in $F(B_i)$ with
isomorphisms
$$\phi_{ij}:X_{j|B_i \times_B B_j} \rightarrow
X_{i|B_i \times_B B_j}$$ in $F(B_i \times_B B_j)$
satisfying the cocycle condition. Then
there is an object $X \in F(X)$ with isomorphisms
$X_{|B_i} \stackrel{\simeq}\rightarrow X_i$ inducing
the isomorphisms $\phi_{ij}$ above.
\end{definition}

Note: If $F$ is a functor, then condition (1) is satisfied
since $Iso_B(X,Y)$ will always be either the empty sheaf
or the constant sheaf. In this case condition (2)
just asserts that the functor is a sheaf in the \'etale
topology.
Condition (2) is not immediate and may easily
fail if $F$ is not representable (A representable functor
will be a stack, since condition
(2) is equivalent to saying that the functor of points
is a sheaf in the \'etale topology). For
example, the moduli functor we defined in Part 1
is not a stack, since it doesn't satisfy condition
(2). In particular, as noted above, given a curve $C$ with
automorphism group $G$ and
$B' \rightarrow B$ Galois with group $G$,  there are two ways to descend
the family $C \times B'/B'$ to a family over $B$,
so a section of $F$ over $B$, is not determined by its pull-back to an \'etale
cover.

However, the moduli groupoid defined above is a stack.
We will not prove this here, but instead
we will prove that the moduli groupoid is
the quotient groupoid of a scheme by $\PGL(N+1)$.

\begin{prop} (cf. \cite[Example 7.17]{Vi})
The groupoid $[X/G]$ defined above is a stack.
\end{prop}
Proof: Let $e,e'$ be sections of $[X/G](B)$
corresponding to principal bundles
$E \rightarrow B$ and $E' \rightarrow B$
and $G$-maps $f: E \rightarrow X$ and $f' : E' \rightarrow X$. Then
$Iso_B(e,e')$ is empty unless $E = E'$ and $f=f'$.
If $e=e'$, then the isomorphisms correspond to elements
of $g$ which preserve $f$. In other words,
$Iso_{B}(e,e)$ is just the $G$-subgroup which
is the stabilizer of the $G-$map $f:E \rightarrow X$
(see \cite[Definition 0.4]{GIT} for the definition of
stabilizer).

The functor which associates to $G$-map $f:E \rightarrow X$
its stabilizer,is represented by the scheme $\mbox{\bf Stab}_X(G)$;
i.e., the stabilizer of the identity map $X \rightarrow X$.
Since $Iso_B(e,e')$ is represented by a scheme it is
a sheaf in the \'etale topology.

Furthermore, a principal $E \rightarrow B$
is determined by \'etale descent, so condition
(2) is satisfied.
\endproof

\subsection{Representable morphisms}

Most of the material in this section is taken
from \cite[Section 4]{DM}.

Let $F$ and $G$ be stacks over $S$. A morphism
of stacks is just a functor of groupoids which commutes
with the projection functor to $\S$.
If $f:F \rightarrow G$ and $h:H \rightarrow G$ are morphisms
of stacks, then we can define the fiber product
$F \times_G H$ as the groupoid whose sections
over a base $B$ are pairs $(x,y) \in F(B) \times H(B)$
such that $f(x)$ is {\it isomorphic} to $h(y)$.
It can be easily checked that this
groupoid is a stack.

\begin{definition}
A morphism $f:F \rightarrow G$ of stacks is
said to be representable if for any map of a scheme
$B \rightarrow G$ the fiber product
$F \times_G X$ is represented by a scheme.
\end{definition}
{\bf Remark:} When we say that a stack is a scheme,
we mean that the stack is the stack associated
to the functor of points of the scheme.

Example: There is a projection map $X \rightarrow [X/G]$
corresponding to the trivial $G$-bundle $X \times_S G$ on $X$.
This morphism is representable because giving a  map $B \rightarrow [X/G]$
is equivalent to giving a principal bundle $E \rightarrow B$.
The fiber product $B \times_{[X/G]} X$ is just the
the scheme $E$.

Let $\P$ be a
property of morphisms of schemes which is stable under
base change and of a local nature on the target.
\begin{definition} (\cite[Definition 4.3]{DM}
A representable morphism of stacks $f:F \rightarrow G$ has property $\P$, if
for all maps of scheme $B \rightarrow G$
the corresponding morphism of schemes $F \times_G B \rightarrow B$
has property $\P$.
\end{definition}
Example: The projection morphism $X \rightarrow [X/G]$ is smooth
since for any $B \rightarrow [X/G]$ the corresponding
map $E \rightarrow B$ is smooth because $E$ is
a principal bundle over $E$.

\subsection{Definition of a Deligne-Mumford stack}
\begin{definition}
A stack is Deligne-Mumford if

(1) The diagonal $\Delta_X:F \rightarrow F \times_S F$ is
representable, quasi-compact and separated.

(2) There is a scheme $U$ and an \'etale surjective
morphism $U \rightarrow F$. Such a morphism
$U \rightarrow F$
is called an (\'etale) atlas.
\end{definition}

{\bf Remark.} In \cite{DM}, such a stack is called
an {\em algebraic} stack. To conform to
current terminology we use the term Deligne-Mumford stack.
A more general class of stacks was studied by Artin,
and they are now called Artin stacks. The basic
difference is that an Artin stack need only have a
{\em smooth} atlas.
Condition (1) above is equivalent to
the following condition:\\

{\it (1') Every morphism $B \rightarrow F$ from a scheme
is representable, so condition (2) makes sense.}

(This fact is stated in \cite{DM} and proved in
\cite[Prop 7.13]{Vi}).

{\bf Remark:} Condition (2) asserts the existence of a
universal deformation space for deformations over Artin
rings.

Vistoli also proves the following proposition:
\begin{prop} \cite[Prop 7.15]{Vi}
The diagonal of a Deligne-Mumford stack
is unramified
\end{prop}

As a consequence of this proposition we can
prove \cite[p. 666]{Vi}
\begin{cor}
If $F$ is a Deligne-Mumford stack, $B$ quasi-compact,
and $X\in F(B)$ then $X$ has only finitely
many automorphisms.
\end{cor}

{\bf Remark.} The are Artin stacks which are not Deligne-Mumford
where each object has a finite
automorphism group. In this case the diagonal is quasi-finite but
ramified.
Objects in the groupoid
have {\em infinitesimal automorphisms.} This phenomenon only occurs
in characteristic $p$, because all groups are smooth in characteristic
0.

Proof: Let $B \rightarrow F$ be map corresponding to
$X$, and let $B \rightarrow F\times_S F$
be the composition with diagonal.
The pullback $B \times_{F \times_S F} F$
can be identified with scheme
$Iso_B(X,X)$. Since $F$ is a Deligne-Mumford stack
the map $Iso_B(X,X)$ is unramified over $X$.
Furthermore, since $B$ is quasi-compact,
the map $Iso_B(X,X) \rightarrow X$ can
have only finitely many sections. Therefore,
$X$ has only finitely many automorphisms over $B$.
\endproof

{\bf Note:} The above proof shows that
diagonal $F \rightarrow F \times_S F$ is not
in general an embedding, since $Iso_B(X,X)$
need not be isomorphic to $B$. It is however,
a local embedding. This is the main technical difficulty
in doing intersection theory directly on Deligne-Mumford
stacks (\cite{Vi}). However, using the equivariant intersection
theory developed in \cite{E-G}, one can avoid these difficulties
for quotient stacks.

The following theorem is stated (but not proved)
in \cite[Theorem 4.21]{DM}. We give the proof
below with a slight additional assumption.
This is the only proof in these notes
which does not appear in the literature.
\begin{thm} \label{smcr}
Let $F$ be a quasi-separated stack over a Noetherian
scheme $S$. Assume that

(1) The diagonal is representable and unramified,

(2) There exists a scheme $U$ of finite type over
$S$ and a smooth surjective $S$-morphism
$U \rightarrow F$,

Then $F$ is a Deligne-Mumford stack.
\end{thm}
{\bf Remark:} This theorem says that condition (1) and
the existence of a versal
deformation space (condition (2)) is actually equivalent to the existence of
a universal deformation space.

{\bf Remark:} We give the proof below under the additional
assumption that the residue fields of the closed points
of $S$ are perfect. In particular we prove the theorem
for stacks of finite type over $Spec\;{\mathbb Z}$.
Using the theorem we will prove that the stack of stable
curves is a Deligne-Mumford stack of finite
type over $Spec\;{\mathbb Z}$.

Proof. The only thing to prove is that $ F$ has an \'etale
atlas of finite type over $S$.
Let $u \in U$ be any closed point in $f^{-1}(u)$. Set
$I_u = \delta^{-1}(u \times_{S} u)$.  Let $z \in U_u$
be a closed point which is separable (i.e. \'etale)
over $u$ (The set of such closed points is dense in a smooth
variety). Since $U_u$ is smooth, the point
$z$ is cut out by a regular sequence in the local ring of $U_u$
at $z$.

The diagonal $\delta: F \rightarrow F \times_S F$
is unramified. Thus, the map $U_u \rightarrow u \times_S U$
obtained by pulling back the morphisms $u \times_S U \rightarrow F \times_S F$
along
the diagonal is unramified. We assume $U$ is of finite type
and that the residue
fields of $S$ are perfect. Thus, $k(u)$ is a finite, hence separable,
extension of the residue field of its image in $S$. Hence
the morphism $u \times_S U \rightarrow U$ is unramified
and so is the composition $U_u \rightarrow u \times_S U \rightarrow U$.

Let $x$ be the image of $z$ in $U$.
By [EGA4 18.4.8] there are \'etale neighborhoods $W'$
and $U'$ of $x$ and $z$ respectively and a closed immersion
$W' \hookrightarrow U'$
such that the diagram commutes
$$\begin{array}{ccc}
W' & \hookrightarrow & U\\
\mbox{{\tiny \'etale}}\downarrow & & \downarrow \mbox{{\tiny \'etale}}\\
U_0 & \rightarrow & U
\end{array}$$

Let $z'$ be any point lying over $y$.
Let $Z_u$ be the closed subscheme of $U'$ defined by lifts to ${\mathcal O}_U'$
of the local equations for $z' \in W'$. By construction, $Z_u$ intersects
$U'$ transversally at $z'$. We will show that the induced
morphism $Z \rightarrow F$ is \'etale in a neighborhood of $Z$.

By definition, this means that for every map of a scheme
$B \rightarrow F$, the induced map of schemes $B \times_F Z_u \rightarrow
B$ is \'etale in a neighborhood of $z' \times_F Z_u$. Since
$U \rightarrow F$ is smooth and surjective, it it suffices to
check that the morphism is \'etale after base change to $U$.

By construction, $Z_u \subset U'$ is cut out by a regular sequence
in a neighborhood of $z' \in U'$ (since $z'$ is a smooth
point of $W'$). Thus $Z_u \times_F U \rightarrow U' \times_F
U$ is a regular embedding in a neighborhood of $z' \times_F U$. Since
$U' \times_F U \rightarrow U'$ is smooth, we can apply
[EGA4,Theorem 17.12.1], and conclude that $Z_u \times_F U' \rightarrow
U'$ is smooth in a neighborhood of $z'$.
Moreover, the relative dimension of this morphism
is $0$. Therefore, $Z_u \rightarrow F$ is \'etale in a neighborhood
of $z'$.

Since $U$ is of finite type over $S$, the $Z_u$'s are as well.
The union of the $Z_u$'s cover $F$ (since their pullbacks via
the morphism $U \rightarrow F$
cover $U$).
Also, $U$ is Noetherian because it is of finite type over
a Noetherian scheme. Thus a finite number of the $Z_u$'s
will cover the $F$. (To see this, we can pullback via the map
$U \rightarrow F$. The pullback of the $Z_u$'s form an \'etale cover
of $U$ which is Noetherian.)
\endproof

\medskip

{\bf Remark:} Without the assumption that the residue fields
of $S$ are perfect we do not know that any point $u \in U$ is
actually unramified over $F$. In this case we need to analyze
the image of $u$ in $F$, which is not a point, but rather
a {\it gerbe} over a point.

\medskip

The theorem has a useful corollary.
\begin{cor}
Let $X/S$ be a Noetherian scheme of finite type
and let $G/S$ be a smooth group scheme acting on
$X$ with finite, reduced stabilizers, then
$[X/G]$ is a Deligne-Mumford stack.
\end{cor}
Proof: The condition on the action
ensures that $Iso_B(E,E)$ is
unramified over $E$ for any map
$B \rightarrow [X/G]$ corresponding
to the principal bundle $E \rightarrow B$.
This in turn implies that the diagonal is
also unramified, so condition (1) is satisfied.
Furthermore, condition (2) is satisfied
by the smooth map $X \rightarrow [X/G]$.
\endproof

\subsection{Further properties of Deligne-Mumford stacks}
Not all morphisms of stacks are representable,
so we can not define algebro-geometric properties
of morphisms as we did for representable morphisms.
However, if we consider morphisms of {\em Deligne-Mumford}
stacks then we can define properties
of morphisms as follows (see \cite[p. 100]{DM})

Let $\P$ be a property of morphisms of
schemes which at source and target is of a local
nature for the \'etale topology. This
means that for any family of commutative
squares
$$\begin{array}{ccc}
X_i & \stackrel{g_i} \rightarrow & X\\
f_i \downarrow & & f \downarrow\\
Y_i & \stackrel{h_i} \rightarrow & Y
\end{array}$$
where the $g_i$ (resp $h_i$) are \'etale and cover $X$ (resp. $Y$),
then $f$ has property $\P$ if and only if $f_i$ has property $\P$
for all $i$.

Examples of such properties are $f$ flat, smooth, \'etale, unramified,
locally of finite type, locally of finite presentation, etc.

Then if $f:F \rightarrow G$ is any morphism of Deligne-Mumford
stacks we say that $f$ has property $\P$ if there are \'etale atlases
$U \rightarrow F$, $U' \rightarrow G$ and a compatible
morphism $U \rightarrow U'$
with property $\P$.

Likewise, if $\P$ is property of schemes which is
local in the \'etale topology (for example regular,
normal, locally Noetherian, of characteristic $p$,
reduced, separated, Cohen-Macaulay, etc.) then a Deligne-Mumford stack
$F$ has property $\P$ if for one (and hence every)
\'etale atlas $U \rightarrow F$, the scheme $U$
has property $\P$.

An open substack $F \subset G$ is a full subcategory of $G$
such that for any $x \in Obj(F)$, all objects in $G$ isomorphic
to $x$ are also in $F$. Furthermore, the inclusion morphism
$F \rightarrow G$
is represented by an open immersions
In a similar way we can talk about closed (or locally closed) substacks.

Using these notions, we say that a map of Deligne-Mumford
stacks $F \rightarrow G$ is separated if for
any map of a scheme $B \rightarrow G$ the fiber
product $F \times_G B$ is separated as a stack over $B$.
It is proper if it is separated,
of finite type and locally over $F$ there is a Deligne-Mumford
stack $H \rightarrow F$ and a representable proper map
$H \rightarrow G$ commuting with the projection to $F$
and the original map $F \rightarrow G$.
$$\begin{array}{ccc}
H &  & \\
\downarrow & \searrow & \\
F & \rightarrow & G
\end{array}$$

{\bf Remark:}
By a theorem of Vistoli \cite[Prop. 2,6]{Vi} and Laumon--Moret-Baily
\cite[Theorem 10.1]{L-MB}
every Noetherian stack has a finite cover by a scheme. Using this
fact we can say that a morphism $F \rightarrow G$ is proper
if there is a finite cover $X \rightarrow F$ by a scheme such
that the composition $X \rightarrow F \rightarrow G$ is a proper
representable morphism. (Recall that any morphism from a scheme
to a stack is representable). Similarly we
say that a morphism $f: F \rightarrow G$ of Noetherian stacks
is (quasi)-finite if for
any finite cover $X \rightarrow F$, the composition
$X \rightarrow F \rightarrow G$ is representable and (quasi)-finite.

As is the case with schemes, there are
valuative criteria for separation and properness
(\cite[Theorem 4.18-4.19]{DM}). The valuative
criterion for separation is equivalent to the
criterion for schemes, but we only construct an {\it isomorphism}
between two extensions.
\begin{thm}\label{sep}
A morphism $f:F \rightarrow G$ is separated iff
the following condition holds:

For
any complete discrete valuation ring $V$ and fraction
field $K$ and any morphism $f:Spec\;V \rightarrow G$ with
lifts $g_1,g_2:Spec\; V \rightarrow F$ which are isomorphic
when restricted to $Spec\;K$, then the isomorphism can be extended
to an isomorphism between $g_1$ and $g_2$.
\end{thm}

\begin{thm}\label{val}
A separated morphism $f: F \rightarrow G$ is proper if and only if
for any complete discrete valuation ring $V$ with field of
fractions $K$ and any map $Spec\;V \rightarrow G$ which lifts
over $Spec\;K$ to a map to $F$, there is a finite separable extension
$K'$ of $K$ such that the lift extends to all
of $Spec\;V'$ where $V'$ is the integral closure of $V$ in $K'$
\end{thm}

{\bf Remark:}
When applied to schemes, Theorem \ref{val} appears to be
stronger than the usual valuative criterion for properness. However,
this is not the case, as it is easy to show that if there
is a lift $Spec\;V' \rightarrow F$, then there is in
fact a lift $Spec\;V \rightarrow F$, as long as
the image of $V'$ is contained in an affine subscheme of $F$ -
which is always the case if $F$ is a scheme.

Finally we conclude this section with the definition of
the moduli space of a Deligne-Mumford stack. This definition
is completely analogous to Mumford's definition (\cite[p. 99]{GIT})
of
a coarse moduli scheme mentioned above.
\begin{definition}
The moduli space of a Deligne-Mumford stack $F$ is a scheme
$M$ together with a proper morphism $\pi:F \rightarrow M$,
such that

(*)  for any algebraically closed field $k$
there is a bijection between the connected components
of the groupoid $F(\Omega)$ and $M(\Omega)$, where
$\Omega = Spec\;k$.

Furthermore, $M$ is universal in the sense that
if $F \rightarrow N$ is a proper map satisfying (*)
there is a morphism $M \rightarrow N$ such that
the map $F \rightarrow N$ factors through $\pi$.
\end{definition}

{\bf Remark:} The reader at this point may wonder
why we need the valuative criterion
for stacks as stated in Theorem \ref{val}.
The difference
is explained as follows. Let $F$ be a complete Deligne-Mumford stack whose
sections are schemes.
If $M$ be a moduli space $F$, then $M$ is also complete.
Let $B= Spec\;{\mathcal O}$ where ${\mathcal O}$ is a DVR with
function field $K$ and residue field $k$, and suppose
there is a map $Spec\; K \rightarrow F$ corresponding to
a section $X_K \rightarrow Spec\;K$ of $F$ over $Spec \; K$.
Since $M$ is complete, the induced map $Spec \; K \rightarrow
M$ can be extended to a map $B \rightarrow M$. However, there need not
be a section $X \rightarrow B$ which restricts to
$X_K\ \rightarrow Spec\; K$. We can only assert that there is
a finite cover $B' \rightarrow B$ and a section $X' \rightarrow B'$
which restricts over the generic fiber to $X_K \times_B B'$.
\section{Stable curves}
In this section we discuss stable curves and the compactification
of the moduli of curves to the moduli of stable curves.
\begin{definition} \cite[Definition 1.1]{DM}
A Deligne-Mumford stable (resp. semi-stable) curve of genus $g$ over a scheme
$S$ is a proper
flat family $C \rightarrow S$ whose geometric fibers are reduced,
connected, 1-dimensional schemes $C_s$ such that:

(1) $C_s$ has only ordinary double points as singularities.

(2) If $E$ is a non-singular rational component of $C$, then
$E$ meets the other components of $C_s$ in more than 2 points (resp.
in at least 2 points).

(3) $C_s$ has arithmetic genus $g$; i.e. $dim\;H^1({\mathcal O}_{C_s}) = g$.
\end{definition}
{\bf Remark:} Clearly, a smooth curve of genus $g$ is stable. Condition
(2) ensures that stable curves have finite automorphism groups,
so that we will be able to form a Deligne-Mumford stack
out of the category of stable curves. We will not use
the notion of semi-stable curves until we discuss geometric
invariant theory in Section 4.

Denote by $F_{\Mgbar}$ the groupoid over $Spec\; \Z$
whose sections over a scheme $B$ are families
of stable curves $X \rightarrow B$. As is the case
with smooth curves, we define a morphism from $X'\rightarrow B'$
to $X \rightarrow B$ as a commutative diagram
$$\begin{array}{ccc}
X' & \rightarrow & X\\
\downarrow & & \downarrow\\
B' &  \rightarrow & B
\end{array}$$
which induces an isomorphism $X' \simeq X \times_B B'$.

\subsection{The stack of stable curves is a Deligne-Mumford stack}

Let $\pi:C \rightarrow S$ be a stable curve. Since $\pi$ is flat
and its geometric fibers are local complete intersections, the
morphism is a local complete intersection morphism. It follows
from the theory of duality that there is a canonical invertible dualizing
sheaf $\omega_{C/S}$ on $C$. If $C/S$ is smooth, then this sheaf
is the relative cotangent bundle. The key fact we need about
this sheaf is a theorem of Deligne and Mumford \cite[p. 78]{DM}.
\begin{thm}
Let $C \stackrel{\pi} \rightarrow S$ be a stable curve of genus $g \geq 2$.
Then $\omega_{C/S}^{\otimes n}$ is relatively very ample for
$n \geq 3$, and $\pi_*(\omega_{C/S}^{\otimes n})$ is locally
free of rank $(2n-1)(g-1)$.
\end{thm}
{\bf Remark:} When $\pi$ is smooth, the theorem follows from
the classical Riemann-Roch theorem for curves. The general
case is proved by analyzing the locally free sheaf
obtained by restricting $\omega_{C/S}$ to the geometric
fibers of $C/S$. In particular, if $S = Spec\; k$, with
$k$ algebraically closed, then $\omega_{C/S}$ can be described as
follows. Let $f:C' \rightarrow C$ be the normalization of $C$
(note $C'$ need not be connected). Let $x_1, \ldots , x_n,y_1, \ldots
, y_n$ be the points of $C'$ such that the $z_i = f(x_i)=f(y_i)$
are the double points of $C$. Then $\omega_{C/S}$ can be identified
with the sheaf
of 1-forms $\eta$ on $C'$ regular except for simple poles at the $x$'s and
$y$'s and with $Res_{x_i}(\eta) + Res_{y_i}(\eta)=0$.

As a result, every stable curve can be realized as a
curve in $\P^{N=(2n-1)(g-1) -1}$ with Hilbert polynomial
$P_{g,n}(t) = (2nt-1)(g-1)$. There is a subscheme (defined over
$Spec\;\Z$)
$\Hbar_{g,n} \subset Hilb^{P_{g,n}}_{\P^N}$ of the Hilbert scheme
corresponding to $n-$canonically embedded stable curves.
Likewise there is a subscheme $H_{g,n} \subset \Hbar_{g,n}$
corresponding to $n-$canonically embedded smooth curves.
A map $S \rightarrow \Hbar_{g,n}$ corresponds to a stable
curve $C \stackrel{\pi} \rightarrow S$ of genus $g$ and an isomorphism
of $\pi_*(\omega_{C/S})$ with $\P^N \times S$.

Now, $\PGL(N+1)$ naturally acts on $H_{g,n}$ and $\Hbar_{g,n}$.
\begin{thm}
$F_{\Mg} = [H_{g,n}/\PGL(N+1)]$ and $F_{\Mgbar}=[\Hbar_{g,n}/\PGL(N+1)]$.
\end{thm}
Note that the theorem asserts that the quotient is independent of $n$.

Proof: Given a family of stable curves $C \stackrel{\pi}\rightarrow B$,
let $E \rightarrow B$ be the principal $\PGL(N+1)$ bundle associated
to the projective bundle $\P(\pi_*(\omega_{C/B}^{\otimes n}))$. Let
$\pi':C \times_B E \rightarrow E$ be the pullback family. The
pullback of this projective bundle to $E$ is trivial and is isomorphic
to $\P(\pi'_*(\omega_{C \times_B E/E}))$, so there
is a map $E \rightarrow \Hbar_{g,n}$ which is clearly $\PGL(N+1)$
invariant. Thus there is a functor $F_{\Mgbar} \rightarrow [\Hbar_{g,n}/
\PGL(N+1)]$ which takes $F_{\Mg}$ to $[H_{g,n}/\PGL(N+1)]$.

The next step is to show that if $C/B$ is a stable curve
then any automorphism of $C/B$ is induced by an
automorphism of the projective bundle $\P(\pi_*(\omega_{C/B}))$.
This is proved for smooth curves in \cite[Proposition 5.2]{GIT},
and is easily
generalized to stable curves because $\pi_*(\omega_{C/B})$
has the same properties as in the smooth case.
It then follows that
$F_{\Mgbar}$ is a full subcategory of the quotient $[H_{g,n}/\PGL(N+1)]$.

Now if $E \rightarrow B$ is a section of $[\Hbar_{g,n}/\PGL(N+1)]$ then
there is a family $C_E \stackrel{\pi_E} \rightarrow E$ of curves of genus $g$
together with an isomorphism $\P(\pi_{E,*}(\omega_{C_E/E})) \simeq \P^N_E$.
Now $\PGL(N+1)$ acts by changing the isomorphism. In particular,
if $g \in PGL(N+1)$ and $e \in E$ then the fiber of $\pi_E$ is
the same over $e$ as it is over $ge$. We can therefore form
a quotient $C/B$ such that $C_E \simeq C \times_B E$.
Hence we have defined a section of $F_{\Mgbar}$,
so there is an equivalence of categories between $F_{\Mgbar}$ and
$[\Hbar_{g,n}/\PGL(N+1)]$ as desired. \endproof

\begin{cor}
$F_{\Mg}$ and $F_{\Mgbar}$ are Deligne-Mumford stacks.
\end{cor}
Proof: We have just shown that $F_{\Mg}$ and $F_{\Mgbar}$
are quotients of a scheme by a smooth group,
so they have smooth atlases. Every stable curve
defined over an algebraically closed field has a finite and reduced
automorphism group, so the diagonal is unramified. Therefore,
they are Deligne-Mumford stacks by theorem \ref{smcr}.
\endproof

\subsection{Properness of $F_{\Mgbar}$}

Given that $F_{\Mgbar}$ is a Deligne-Mumford stack, the valuative
criterion of properness and the following stable reduction
theorem show that it is proper over $Spec\;\Z$.
\begin{thm} \label{stable}
Let $B$ be the the spectrum of a DVR with function field
$K$, and let $X \rightarrow B$ be a family of curves
such that its restriction $X_K \rightarrow Spec\; K$ is
a stable curve. Then there is a finite extension $K'/K$
and a unique stable family $X' \rightarrow B \times_K K'$
such that the restriction $X' \rightarrow Spec\; K'$
is isomorphic to $X_K \times_K K'$.
\end{thm}
Remarks on the proof of Theorem \ref{stable}: This theorem
was originally proved (but not published)
in characteristic zero by Mumford and
Mayer (\cite[Appendix D]{GIT}). There is a relatively straightforward
algorithmic version of this theorem in characteristic 0 which I
learned from Joe Harris.
Blow up the singular points of the special fiber of $X/B$
until the total space of the family is smooth
and the special fiber has only nodes as singularities. The modified
special fiber will have a number of components with
positive multiplicity coming from the exceptional divisors
in the blowups. Next, do a base change of degree
equal to the g.c.d. of the multiple components. After base change
all components of the special fiber will have multiplicity 1.
Then contract all (-1) and (-2) rational components in
the total space. The special fiber is now stable. Furthermore,
the total space of the new family is a minimal model for the
surface. Since minimal models of surfaces are unique,
the stable limit curve is unique. \endproof

This algorithmic proof fails in characteristic $p>0$, because
after blowing up some components of the special fiber may
have multiplicity divisible by $p$. In this case,
it will not be possible to make the component
become reduced after base change.

Deligne and Mumford proved the stable reduction theorem in arbitrary
characteristic using Neron models of the Jacobians of the
curves (\cite{DM}). Later Artin and Winters
gave a direct geometric proof using the theory of curves
on surfaces.

\subsection{Irreducibility of $F_{\Mg}$ and $F_{\Mgbar}$}

Using the description of the moduli stacks as quotients of
$H_g$ and $\Hbar_g$ we can deduce properties of the stacks
from the corresponding properties of the Hilbert scheme.
In particular, deformation theory shows that $H_g$ and
$\Hbar_g$ are smooth over $Spec\;\Z$ (\cite[Cor 1.7]{DM}).
Since the map $H_g \rightarrow F_{\Mg}$ (resp.
$\Hbar_g \rightarrow F_{\Mgbar}$) is smooth we see that
$F_{\Mgbar}$ is smooth.

Further analysis \cite[Cor 1.9]{DM} shows that the scheme
$\Hbar_g - H_g$ representing polarized,
singular, stable curves is a divisor with
normal crossings in $\Hbar$. This property descends to
the moduli stacks.

\begin{thm} \cite[Thm 5.2]{DM} $F_{\Mgbar}$ is smooth and proper
over $Spec\; \Z$. The complement
$F_{\Mgbar} - F_{\Mg}$ is a divisor with normal crossings
in $F_{\Mgbar}$.
\end{thm}

The main result of \cite{DM} is the following theorem:
\begin{thm} \cite{DM} \label{irred}
$F_{\Mgbar}$ is irreducible over $Spec\; \Z$.
\end{thm}

{\bf Remark:}  Deligne and Mumford gave two proofs of this theorem. In both
cases they deduce the result from the classical
characteristic 0 result stated below.
We outline below their second proof, which uses
Deligne-Mumford stacks.

\begin{prop}
$F_{\Mg} \times_{Spec\;\Z} Spec\;\C$ is irreducible.
\end{prop}
Proof: It was shown classically that there is a space $H_{k,b}$
parametrizing degree $k$ covers of $\P^1$ simply branched over $b$ points
defined over the complex numbers. In
\cite{Fu}, Fulton showed that
the functor $F_{H_{k,b}}$ whose sections over a base $B$
are families of smooth curves $C \rightarrow B$ together with a
degree $k$ map $C \rightarrow \P^1_B$ expressing each geometric fiber
as a cover of $\P^1$ simply branched over $b$ points
is represented
by a scheme which we also call $H_{k,b}$.
In characteristic greater than $k$ it  is a
finite \'etale cover of $P_b = (\P^1)^b - \Delta$, where
$\Delta$ is the union of all diagonals (This fact was known
classically over $\C$. In low characteristic the map
may fail to be finite). Since $P_b$ is obviously
irreducible, it can be proved that $H_{k,b}$ is irreducible in
high characteristic
by showing that the monodromy of the covering $H_{k,b} \rightarrow
P_b$ acts transitively on the fiber over a base point in $P_b$ for
all $k,b$.
Since there is a universal family of branched covers
$C_{k,b} \rightarrow H_{k,b}$ there is a map $H_{k,b} \rightarrow
F_{\Mg}$ (where $g=b/2 -k+1$). By the Riemann-Roch theorem for
smooth curves, every curve of genus $g$ can be expressed
as a degree $k$ cover of $\P^1$ with $b$ simple branch points,
as long as $k>g+1$. Thus for $k$ (and thus $b$) sufficiently large,
the map is surjective. Therefore $F_{\Mg}$ is irreducible in characteristic
greater than $k$, and thus $F_{\Mg} \times_{\Z} \C$ is irreducible.
\endproof

Proof of Theorem \ref{irred}(outline):
Since $F_{\Mgbar} - F_{\Mg}$ is a divisor
$F_{\Mgbar}$ is irreducible if and only if
$F_{\Mg}$. The stack  $F_{\Mg}$ is smooth, so it suffices to
show that it is connected.

In \cite[Section 5]{DM}, Deligne and Mumford construct
a stack $_{n}F_{\Mg}$ whose sections are curves with a ``Jacobi structure
of level $n$ '' (\cite[Paragraph 5.14]{DM}). The extra structure
eliminates all non-trivial automorphisms when $n \geq 3$, so this stack is in
fact represented by a scheme. Now, $_{n}F_{\Mg}$ is a finite cover
of $F_{\Mg}$, and they deduce the connectedness of $_{n}F_{\Mg} \times \C$
from the connectedness of $F_{\Mg} \times \C$ (\cite[Theorem 5.13, Lemma 5.16]
{DM}). They also prove that the fibers of $_{n}F_{\Mg} \rightarrow
Spec\;\Z[e^{2\pi i/ n},1/n]$ have the same number of connected
components \cite[Cor 5.11]{DM}.
Since $n \geq 3$ was arbitrary, the irreducibility
of $F_{\Mg}$ over $Spec\;\Z$ follows.
\endproof

{\bf Remark:} In \cite{HM} Harris and Mumford constructed a compactification
of $H_{k,b}$ where the boundary represents stable curves expressed as
branched covers of chains of $\P^1$'s. The existence of this
compactification implies that every smooth curve admits degenerations
to singular stable curves. Fulton \cite{Fu82} used this fact to resurrect an
argument of Severi
to give a purely algebraic proof that $F_{\Mg}$ is irreducible
in characteristic 0. This combined with the results of
\cite{DM} give a purely algebraic proof
that $F_{\Mg}$ is irreducible over $Spec\;\Z$.

\section{Construction of the moduli scheme}
As we have previously seen, the moduli stack is a quotient
stack of a smooth scheme $H_g$ by $\PGL(N+1)$. In this, the final
section, we discuss the construction of a quotient scheme
$H_g/\PGL(N+1)$ over an algebraically closed
field $k$. We first prove that such a scheme is unique
and is the coarse moduli space for the quotient stack. We then
briefly discuss Gieseker's GIT construction of a quotient scheme.

Throughout this section, we will assume that all schemes
are defined over an algebraically closed field $k$.
\subsection{Geometric quotients}
\begin{definition} \cite[Definitions 0.5, 0.6]{GIT}
Let $X$ be a scheme defined over a field $k$, and let
$G/k$ be an algebraic group  acting on
$X$. A $k$-scheme $Y$ is a {\em geometric quotient} of
$X$ by $G$ if there is a morphism $X \rightarrow Y$ such
that

(1) f is surjective, affine and $G$ invariant.

(2) $f_*({\mathcal O}_X)^G = {\mathcal O}_Y.$

(3) If $W \subset X$ is closed and $G$ invariant, then $f(W)$
is closed in $Y$. Furthermore, if $W_1$ and $W_2$ are disjoint
$G$ invariant subsets of $X$, then $f(W_1)$ and $f(W_2)$
are disjoint.

(4) The geometric fibers of $f$ are orbits.
\end{definition}

{\bf Remark:} The purpose of the geometric invariant theory
developed by Mumford
is to construct geometric quotients for the action of
reductive\footnote{The definition of a reductive
group is given in \cite[Appendix A]{GIT}. For the purpose
of these notes, it suffices to know that
$SL(N+1,k)$ is reductive for a
field $k$.} groups on projective varieties.

The following is a restatement of \cite[Prop 0.1]{GIT}.
\begin{prop} \label{cat}
If $X  \stackrel{f} \rightarrow Y$ is a geometric quotient and if
$X \stackrel{g} \rightarrow Z$ is a $G$ invariant morphism, then
there is a unique morphism $\phi:Y \rightarrow Z$ such that
$g=\phi \circ f$.
\end{prop}
Note that the proposition implies that if a geometric
quotient exists then it is unique.

Now let $X$ be a scheme with a $G$ action such that the stabilizers of
geometric points are finite and reduced.
We have seen that the groupoid $[X/G]$ is a Deligne-Mumford
stack. Assume the action of $G$ on $X$
is locally proper; i.e., $X$ can be covered
by open invariant subschemes $U$ such that $G$ acts properly on
$U$.
\begin{prop} \cite[Proposition 2.11]{Vi}
If $f: X \rightarrow Y$ is a geometric quotient
of $X$ by $G$, and $f$ is universally submersive
\footnote{$f:X \rightarrow Y$ is submersive if $U \subset Y$
is open if and only if $f^{-1}(U)$ is open in $X$. If this
property remains after base change, then we say
$f$ is universally submersive.},
then $Y$ is the (coarse) moduli space for the
stack $[X/G]$.
\end{prop}
Proof: Let $\eta \in [X/G](B)$ be a section corresponding
to a principal $G$-bundle $E \stackrel{\pi} \rightarrow B$ together with
a $G$ invariant map $\phi:E \rightarrow X$. Since the map
$f: X \rightarrow Y$ is $G$ invariant, there is a unique map
$\psi:B \rightarrow Y$ making the obvious square commute.
Hence we can associate to $\eta$ a unique section of
of $Hom(B,Y)$. Thus there is a morphism of stacks
$[X/G]  \rightarrow Y$.

Now if $\Omega = Spec\;K$ where $K$ is algebraically closed,
then $Hom(\Omega, Y)$ is, by Condition (4) of the definition,
the set of orbits of $K$-valued points of $X$. This is exactly
$[X/G](\Omega)$. Therefore, the map induces a bijection
$[X/G](\Omega) \rightarrow Y(\Omega)$ as required in the
definition of a coarse moduli space.

Using the local properness of the $G$ action together with the universal
submersiveness of $f$ we can prove that the morphism
$[X/G] \rightarrow Y$
is proper (\cite[Proposition 2.11]{Vi}.

Finally, note that the universal property of $Y$ implied
by Prop \ref{cat} also shows that $Y$ satisfies
the universal property necessary for a (coarse) moduli space.
\endproof

{\bf Remark:} If $G$ is reductive and $f:X \rightarrow Y$
is a geometric quotient, so that $f$ is affine, then
\cite[Proposition 0.8]{GIT} the action of $G$ is actually
proper. This condition is satisfied in all geometric invariant
theory quotients of quasi-projective varieties.

\subsection{Construction of quotients by geometric invariant
theory}

In this paragraph we discuss the geometric invariant
theory necessary to construct $\Mg$ and
and $\Mgbar$ as quotients of Hilbert schemes of
$n$-canonically embedded (stable) curves. Our source
is \cite[Chapter 0]{Gi}. For a full
treatment of geometric invariant the classic reference is Mumford's
\cite{GIT}.

Let $X \subset \P^N$ be a projective scheme, and let
$G$ be a reductive group
acting on $X$ via a representation $G \rightarrow GL(N+1)$.

\begin{definition}
(1) A closed point $x \in X$ is called semi-stable
if there exists a non-constant $G$-invariant homogeneous polynomial
$F$ such that $F(x) \neq 0$.

(2) $x \in X$ is called stable if: $dim\;o(x) = dim\; G$ (where
$o(x)$ denotes the orbit of $x$) and there exists
a non-constant $G$-invariant polynomial such that $F(x) \neq 0$
and for every $y_0$ in $X_F=\{y \in X | F(y) \neq 0\}$,
$o(y_0)$ is closed in $X_F$.
\end{definition}

Let $X^{ss}$ denote the semi-stable points of $X$, and
$X^{s}$ denote the stable points. Then $X^{s} \subset X^{ss}$
are both open in $X$. However, they may be empty.

The following is the first main theorem of geometric invariant
theory.
\begin{thm}
There exists a projective scheme $Y$ and
a universally submersive
morphism $f_{ss}:X^{ss} \rightarrow Y$ such that
$f_{ss}$ satisfies properties (1)-(3) of the definition
of a geometric quotient (such a morphism is often
called a good quotient in the literature).
Furthermore, there exists $U \subset Y$ open such that
$f^{-1}(U) =X^{s}$ and $f_s:X^s \rightarrow U$ is
a geometric quotient of $X^s$ by $G$.
\end{thm}

\subsection{Criteria for stability}
Let $X \subset \P^N$ be a projective scheme, and let $\tilde{X}
\subset \A^{N+1}$ be the affine cone over $X$. Assume as above,
that a reductive group $G$ acts on $X$ via a representation
$G \rightarrow GL(N+1)$. Then $G$ acts on $\tilde{X}$
as well. The stability of $x \in X$ can be rephrased in
terms of the stability of the points $\tilde{x} \in
\tilde{X}$ lying over $x$.
\begin{prop} \cite[Chapter 1, Proposition 2.2 and Appendix B]{GIT}
A geometric point $x \in X$ is semi-stable if
for one (and thus for all) $\tilde{x} \in X$ lying
over $X$, $0 \notin o(\tilde{x})$.
The point $x$ is stable if $o(\tilde{x})$ is
closed in $\A^{N+1}$ and has dimension equal to the dimension of $G$.
\end{prop}

The second main theorem of geometric invariant theory
is Mumford's numerical criterion for stability which we
now discuss.
\begin{definition}
A 1-parameter subgroup of $G$ is a homomorphism
$\lambda: G_m \rightarrow G$. This will be abbreviated to
$\lambda$ is a 1-PS of $G$.
\end{definition}
Now if $\lambda$ is a 1-PS of $G$, then the since
$\lambda$ is 1-dimensional, there is a basis
$\{e_0, \ldots , e_{N} \}$  of
$\A^{N+1}$ such that the action of $\lambda$ is diagonalizable
with respect to this basis; i.e. $\lambda(t)e_i = t^{r_i}e_i$
where $t \in G_m$ and $r_i \in \Z$. If $\tilde{x}=
\sum x_ie_i \in \tilde{X}$, then the set of
$r_i$ such that $x_i$ is non-zero is called the $\lambda$-weights
of $\tilde{x}$. Note that if $x \in \P^N$ then the
$\lambda$-weights
are the same for all points in $\A^{N+1}-0$ lying over $x$.

\begin{definition}
$x \in X$ is $\lambda$-semi-stable if for one
(and thus for all) $\tilde{x} \in \tilde{X}$ lying
over $x$, $\tilde{x}$ has a non-positive $\lambda$
weight.
A point $x$ is $\lambda$-stable if $\tilde{x}$ has a negative
$\lambda$-weight.
\end{definition}

\begin{thm} \cite{GIT}
A point $x \in X$ is (semi)stable if and only if $x$ is
$\lambda$-(semi)stable for all 1-PS $\lambda: G_m \rightarrow
G$.
\end{thm}
Remark on the Proof: It is easy to see that
if $x$ is unstable (i.e. not semi-stable)
with respect to $\lambda: G_m \rightarrow G$
then $x$ is unstable. The reason is that
if all the weights of $\lambda$ are positive
then 0 will be in the closure of the
$G$-orbit of $\tilde{x}$ in $\A^{N+1}-0$.
The converse is more difficult.
\endproof

Example (cf. \cite[Proposition 4.1]{GIT}).
The set of homogeneous forms of degree 4 in two
variables forms a 5-dimensional vector space $V$.
We will view $\P(V)$ as the space parametrizing
4-tuples of (not necessarily) distinct points
in $\P^1$. There is a natural action of $SL(2)$
on $V$ inducing an action on $\P(V)$. Let us
use the numerical criterion to determine the
stable and semi-stable locus in $\P(V)$.

If $v \in V$ is a form of degree 4 and
$\lambda$ is a 1-PS subgroup of $SL(2)$,
then we can write $v=a_4X_0^4 + a_3X_0^3X_1 + a_2X_0^2X_1^2
+ a_1X_0X_1^3 + a_0X_1^4$, and
$\lambda$ acts by $\lambda(t)(X_0)=t^rX_0$,
$\lambda(t)(X_1) = t^{-r}X_0$  and $r > 0$(the weight
on $X_1$ must be the negative of the weight on
$X_0$, since $\lambda$ maps to $SL(2)$).
The possible weights of $v$ are
$\{4r,2r,0,-2r,-4r\}$. In order for
$v$ to be $\lambda$-stable one of
$a_1$ or $a_0$ must be non-zero. It is $\lambda$-semi-stable
if one of $a_2$, $a_1$ or $a_0$ is non-zero.
On the other hand, we can consider the 1-PS, $\tau$
which acts by $\tau(t)X_0 = t^{-r}X_0$ and
$\tau(t)X_1 = t^rX_1$. In order for $v$ to be
$\tau$-stable one of $a_4$ and $a_3$
must be non-zero, while it is $\tau$-semi-stable
if $a_2$ is non-zero. Combining the
conditions imposed by $\lambda$ and $\tau$ we see
that if $v$ is stable, then one of $a_0$ or $a_1$
is non-zero and one of $a_3$ or $a_4$ is non-zero.
This condition is equivalent to the condition
that $(1:0)$ and $(0:1)$ are not multiple points of
the subscheme of $\P^1$ cut out by the form $v$.
Likewise, $v$ is semi-stable if
$(1:0)$ or $(0:1)$ is cut out with multiplicity
no more than 2. Finally $v$ is unstable if
$(1:0)$ or $(0:1)$ is cut out with multiplicity
more than 2.

{}From this analysis it is clear that if
$v \in V$ cuts out 4 distinct points then it
will be stable for every 1-PS. Likewise if
$v$ cuts out a subscheme of $\P^1$
with each point having multiplicity 2 or less
then it is semi-stable for every 1-PS. Conversely,
if $v$ cuts a point of multiplicity 3 or more
then $v=X_0^3(a_0X_0 + a_1X_1)$ for some choice
of coordinates on $\P^1$. Then $v$
will have strictly positive weights for
a 1-PS $\lambda$ acting diagonally by $\lambda(t)X_0 = t^rX_0$
for $r>0$.

\subsection{Gieseker's construction of $\Mgbar$}
Let $Hilb^{N+1}_{P(t)}$ be the Hilbert scheme of
curves in $\P^N$ with Hilbert Polynomial $P(t)$.
Now if $X \subset \P^N$ is a curve with Hilbert polynomial
$P(t)$, then there exists $m >> 0$ (independent of $X$) such that
the restriction map $H^0(\P^N, {\mathcal O}_{\P^N}(m)) \rightarrow
H^0(X,{\mathcal O}_X(m))$ is surjective and
$\mbox{dim }H^0(X,{\mathcal O}_X(m)) = P(m)$.
Taking the $P(m)$-th exterior power of $\phi_m$ we
obtain a linear map $V^m = \bigwedge^{P(m)} H^0(\P^N, {\mathcal O}_{\P^N}(m))
\rightarrow \bigwedge^{P(m)} H^0(X,{\mathcal O}_X(m)) = k$ unique up to
scalars; i.e.,  an element of $\P(V^m)$. The corresponding point
in $\P(V^m)$ is called the $m$-th Hilbert point of $X$ and
is denote $H_m(X)$. In this
way we obtain a map $Hilb^{N+1}_{P(t)} \rightarrow \P(V^m)$. For
$m$ sufficiently large this map is an embedding.

Both $SL(N+1)$ and $\PGL(N+1)$ act on $\P(V^m)$ via the $m$-th
exterior power representation of $SL(N+1) \rightarrow GL(m)$.
Now the action of $SL(N+1)$ factors through the action of
$\PGL(N+1)$ (the stabilizer of $SL(N+1)$ at a geometric
point is the group of $N+1$ roots of unity) so we have
the following proposition.
\begin{prop} \label{gpfc}
If $X \subset \P(V^m)$ then
$X \rightarrow Y$ is a geometric quotient by
$SL(N+1)$ if and only if it is a geometric quotient
by $\PGL(N+1)$.
\end{prop}
Proof: If $X \rightarrow Y$ is a geometric quotient by $SL(N+1)$
then the geometric fibers are $SL(N+1)$ orbits. These orbits
are the same as the $\PGL(N+1)$ orbits. Likewise, ${\mathcal O}_X^{SL(N+1)}
= {\mathcal O}_X^{\PGL(N+1)}$. Thus, ${\mathcal O}_Y \simeq f_*{\mathcal
O}_X^{\PGL(N+1)}$.
Finally if $W$ and $V$ are $\PGL(N+1)$ invariant, they are also $SL(N+1)$
invariant. Thus if they are disjoint, then since $X \rightarrow Y$ is
an $SL(N+1)$ quotient, their images will be disjoint as well. Hence
$X \rightarrow Y$ is a $\PGL(N+1)$ quotient.
The converse is similar.
\endproof

Let $g \geq 3$ and $d \geq 20(g-1)$ be integers. Consider the Hilbert
scheme $H^{N+1}_{P(t)}$ of curves in $\P^{N=d-g}$ with Hilbert polynomial
$P(t) = dt-g+1$ (the curves parametrized necessarily
have arithmetic genus $g$). The first step in Gieseker's
construction is to prove the following theorem. The proof
is 10 pages long and uses the numerical criterion.
\begin{thm} \cite[Theorem 1.0.0]{Gi}
There exists an integer $m_0 >> 0$ such that if $X$ is smooth then
$H_{m_0}(X)$ is $SL(N+1)$ stable.
\end{thm}
{\bf Remark:} The theorem is not necessarily true for arbitrary
$m_0>> 0$. However there are infinitely many $m_0$
for which the theorem is true (\cite[Remark after Theorem 1.0.0]{Gi}).

The next, and technically most difficult step is to prove
the following theorem. The proof takes 50 pages!
\begin{thm} \cite[Theorem 1.0.1]{Gi}
For the same integer $m_0$, every point in $Hilb^{N+1}_{P(t)}
\cap \P(V^{m_0})^{ss}$ parametrizes a Deligne-Mumford semi-stable
curve.
\end{thm}

Let $U \subset Hilb_{P(t)}^{N+1}$ be the subscheme of
semi-stable curves with respect to the $m_0$-th Hilbert
embedding. Let $Z_U \subset \P^N_U$ be the restriction of the universal
family of projective curves. As before, view a point $h \in U$
as parametrizing a curve $X_h$ and a very ample line bundle
$L_h$ of degree $d$ on $X_h$.
Set $U_c = \{h \in U | L_h \simeq \omega_{X_h}^n\}$. This is a constructible
subscheme of $U$ which is empty unless $2g-2$ divides $d$.
Gieseker then proves that $U_c$ is in fact closed in $U$. He also
proves that $U_c$ is smooth (\cite[Theorem 2.0.1]{Gi})
and parametrizes only all Deligne-Mumford stable curves; thus,
$U_c \simeq \overline{H}_{g,n}$. Since
$U_c$ is closed in $U$ there is a projective quotient $U_c/SL(N+1)$.
Finally note that $\PGL(N+1)$ (and thus $SL(N+1)$) acts with finite
stabilizers on points of $U_c$ because the curves parametrized
have finite automorphism groups. Hence the points of $U_c$ are in
fact $SL(N+1)$ stable. Thus a universally submersive geometric
quotient $U_c/SL(N+1)$ exists. Since this is isomorphic to
a geometric quotient $U_c/PGL(N+1) \simeq \overline{H}_{g,n}/\PGL(N+1)$
we have succeeded in constructing a coarse moduli scheme for
the stack of stable curves.
\endproof

\end{document}